\documentclass[a4paper]{article}

\usepackage{babel}
\usepackage{amsfonts} % amsthm
\usepackage{amssymb,latexsym,amsmath}
\usepackage{layout}
\usepackage{graphics}
\usepackage{graphicx}
\usepackage{calc}
\usepackage{color}
\usepackage{mathrsfs}
\usepackage{url}
\usepackage{hyperref}

\newcommand\s[1]{\text{Sing}(#1)}

\usepackage[symbol]{footmisc}

\begin{document}

\title{Some mathematical and computational relations between timbre and color}

\author{Maria Mannone\footnote{Department of Engineering, University of Palermo, Italy, mariacaterina.mannone@unipa.it; Dipartimento di Scienze Ambientali, Informatica e Statistica (DAIS) and European Centre for Living Technology (ECLT), Ca' Foscari University of Venice, Italy, maria.mannone@unive.it} \, and
Juan Sebastián Arias-Valero\footnote{Independent researcher, jsariasv1@gmail.com}}

\maketitle

\begin{abstract}
In physics, timbre is a complex phenomenon, like color. Musical timbres are given by the superposition of sinusoidal signals, corresponding to longitudinal acoustic waves. Colors are produced by the superposition of transverse electromagnetic waves in the domain of visible light. Regarding human perception, specific timbre variations provoke effects similar to color variations, for example, a rising tension or a relaxation effect. We aim to create a computational framework to modulate timbres and colors. To this end, we consider categorical groupoids, where colors (timbres) are objects and color variations (timbre variations) are morphisms, and functors between them, which are induced by continuous maps. We also sketch some gestural variations of this scheme. Thus, we try to soften the differences and focus on the similarity of structures. 


\end{abstract}

\begin{quote}
\begin{small}
\end{small}
\end{quote}

%\tableofcontents

\section{Introduction}\label{introduction}\label{intro}

Timbres and colors fascinated musicians, artists, and scientists across centuries \cite{leonardo,kandinsky}. In physics, the complexity of timbre is due to the superposition of simple components (sinusoidal waves), which can be separated with Helmholtz resonators \cite{hemholtz}. Timbres can be computationally investigated with Fourier transforms and sonograms, which show the strength of each component (partial) of the superposition. Colors are also related to the idea of superposition, as proved by \cite{newton} for white light, which can be decomposed in colors through a prism. The physics involved is quite different: sound involves mechanical longitudinal waves, while light is made of electromagnetic transverse waves. However, timbres and colors have a main similarity: they are complex signals, made of simple superposed wave signals. This suggests a correspondence based on the relation between sound frequency and spatial frequency of light, but according to \cite{caivano}, absolute correspondences between these domains are difficult to establish, so the relativity and the obstructions of this problem could be softened in the categorical context.

In fact, the point of view of precise measurement can be enriched in several ways. Scholars such as Goethe pointed out the importance of perception to understand colors in the framework of nature and the arts \cite{goethe}. In addition, both colors and timbres can be qualitatively rated as, for example, {\em cold}, {\em strong}, or {\em delicate}. Even though different cultures can associate a different (symbolic) meaning to each color, we can find aspects with certain universality, related to human perception. Some colors are more instinctively associated with higher or lower tension: red or yellow raise more attention than light blue or gray. Similarly, specific orchestral timbres are more awakening than others: a loud\footnote{Loudness in music performance can affect timbre characteristics \cite{castellengo}.} trumpet sound is a more effective alarm than a soft flute melody. Some recent studies point out the importance of a ``shared emotion'' to associate colors and musical sequences \cite{palmer}, which also occurs in the framework of classical music listening \cite{crnjanski}. On the other hand, both colors and timbres can be mixed or shaded---as it happens for painting and orchestration, respectively, transforming a {\em delicate sound (or color)} into a {\em strong sound (or color)}. In this way, we can draw upon the idea of {\em superposition} and {\em similarity of perception} to imagine how we can investigate colors and timbres, focusing on common aspects through abstraction.\footnote{A first experiment, where participants were asked to associate colors, color bands, and timbres, confirmed a non-negligible perceptive correlation \cite{mannone_santini}.} These aspects are intensities, mixing, and shadows/nuances. In particular, harmonic choices, which influence timbre, are also ruled by the idea of superposition.

In this article we introduce fundamental groupoids of color and timbre spaces and functors between them. These functors could be induced by some classical (possibly) continuous maps suggested in \cite{caivano}. This categorical framework \cite{lurie,CatWork} could be adequate to express the superposition and similarity principles to understand the color/timbre relation, complementing analytical approaches. Categories have already been used to investigate processes and phenomena in the arts from a bird's-eye perspective \cite{jed,Kubota,MazzolaTopos}.

This article is structured as follows. In Section \ref{general_}, we review some color spaces, timbre spaces, and maps between them. In Section \ref{categorical_enrichment}, we offer a categorical enrichment of previous approaches to relate color and timbre. In particular, in Section \ref{example}, we include a computational example of interaction between color and timbre paths. Then, Section~\ref{gest} is devoted to a gestural extension of the previous enrichment. In Section \ref{conc} some conclusions and further possible applications are discussed. In the {\em Glossary} (Section~\ref{G}) we provide definitions of some specialized mathematical concepts that we mention. We use \textbf{boldface} for these terms. 
   
%As a general disclaimer, we make it precise that the topic of colors and timbres and their relations is particularly broad.
As a general disclaimer: colors, timbres, and their relationships constitute a vast topic. This is a position paper (or rather, a working one) aiming to open the way toward further studies in this field.

\section{Spaces and mappings: an overview}\label{general_}
\subsection{The CIE 1931 color space}\label{color}

The CIE model \cite{how} connects the visible spectrum with human perception. It assigns to each spectrum wavelength $\lambda\in [380,780]$, measured in nanometres, three sensitivity level values $\overline{x}(\lambda)$, $\overline{y}(\lambda)$, and $\overline{z}(\lambda)$ corresponding to the kinds of human cone cells under certain standard conditions. Thus, a spectral distribution yields, by integration of its product with each color matching function ($\overline{x}$, $\overline{y}$, or $\overline{z}$), a triple $(X,Y,Z)$ of color coordinates. All these triples amount to the unit cube $[0,1]^3$, after normalizing units. We embed this cube in $\mathbb{R}^3$, regarding the latter as a vector space and a topological space. The vector sum in the cube, whenever defined, corresponds to color mixing (superposition of light beams).\footnote{We mix colors in printing and painting with the subtractive model, a sort of dual of the additive one.} If the sum is not in $[0,1]^3$, one can take an average of vector components to represent a mixture (with average intensity) for computational purposes.

On the other hand, the standard RGB space is used for screens and photography, so we need it for experiments. It has red, green, and blue as primary colors, which give white if superposed. The standard RGB model does not cover the CIE gamut in principle, for instance, \textcolor{black}{a} spectral violet. % is not additive and 
However, we can transform CIE to standard RGB by means of an appropriate conversion of CIE to linear RGB followed by electro-optical transfer. The RGB space has already been considered for mathematical modeling \cite{provenzi,resnikoff}. In particular, \cite{resnikoff} proposed a three-dimensional space of perceived colors, where equivalence classes correspond to perceptual match.

\subsection{Timbre space}

As a possible representation of timbres, we can consider the space proposed by Grey \cite{grey}, based on the dissimilarity between pairs of musical instrument sounds. 

On the other hand, we have the set of all continuous periodic maps. These maps represent continuous sound waves that can be recovered from their Fourier series according to Fourier's and Fejér's theorems \cite[Section~2.4]{offering}. It is embedded in the space of all continuous maps $\mathbb{R}\longrightarrow\mathbb{R}$, which has the \textbf{compact-open topology} and is a vector space. Superposition of waves corresponds to addition of the associated periodic functions, although the result need not be periodic. In what follows, we take the topological space of continuous periodic maps as our timbre space, given the structural analogy with the CIE space in the sense that color/wave superpositions correspond to vector sums. 

\subsection{Maps between timbre and color}\label{map}

According to \cite{caivano}, a possible correspondence between color and sound can be based on the idea that a musical octave should match a color octave. A musical octave is a closed interval of the form $[f,2f]$, where $f$ is a fixed sound frequency in Hertz. Human vision barely ranges through color octave, namely the interval of wavelengths in nanometres $[380,760]$, which corresponds to the interval of spatial frequencies $[(1/2)(1/380),1/380]$ by means of the assignment $\lambda\mapsto 1/\lambda$. Thus, the map $\lambda\mapsto 760f/\lambda$ is a continuous bijection from the color octave $[380,760]$ to the musical octave $[f,2f]$. Note that under this logic, the color order violet-blue-green-yellow-orange-red corresponds to a decreasing pitch frequency.

Since human hearing ranges frequencies in the Hertz interval $[20,20000]$, and therefore several octaves, there is not a perfect correspondence between wavelengths and frequencies. This suggests reducing the interval  $[20,20000]$ modulo a chosen octave and then using the previous correspondence. The resulting map is continuous under the assumption that we identify the endpoints of $[380,760]$.  

There are other possibilities for a correspondence between color and sound. Some scholars focus on perceived correspondences of pitch classes with classes of colors \cite{nature}. The use of classes can be formalized by means of quotient spaces. Classes take into account perceptive similarities but not perfect one-to-one associations. Other continuous correspondences could associate the transition from violet to red with an increasing pitch frequency. 

The following construction is a possible way to get a continuous\footnote{We do not have a proof of this continuity.} map from the timbre space to the CIE color space. First, let us consider the case of a timbre given by simple FM synthesis \cite[Section~8.8]{offering}, namely a periodic\footnote{The wave is periodic if the quotient between carrier and modulator frequencies is a rational number.} wave corresponding to
\begin{equation}\label{bessel1}
\sin[\omega_c t+I\sin(\omega_m t)],
\end{equation}
where $\omega_c=2\pi f_c$, $\omega_m=2\pi f_m$, $f_c$ is the carrier frequency, $f_m$ is the modulator frequency, and $I$ is the modulation index. An associated convergent series is
\begin{equation}
\sum_{n =-\infty}^{\infty}J_n(I)\sin[(\omega_c+n\omega_m)t],
\end{equation} 
where $J_n$ is the $n$th Bessel function of the first kind. This series expresses the wave in terms of simple harmonics with frequencies $f_c+n f_m$ for $n\in \mathbb{Z}$.
% By manipulating this series so that all frequencies are positive we obtain a new one of the form
By factorizing the sign of each negative value of $f_c + nf_m$ outside of $\sin[(f_c + nf_m)t]$ we obtain:
\begin{equation}\sum_{n=0}^{\infty}a_n\sin(2\pi f_n t).
\end{equation}
Thus, given a continuous map $h$ from frequencies to color wavelengths, we construct (by linearity) the series in the CIE space
\begin{equation}\sum\limits_{n=0}^{\infty}a_n XYZ(h(f_n)),
\end{equation}
where $XYZ(\lambda)$ gives the CIE coordinates of the wavelength $\lambda$, whenever the series converges in the CIE space. In general, one could use the Fourier series \cite[p.~1019]{MazzolaTopos}:
\begin{equation}\label{timbre_coordinates}
a_0+\sum\limits_{n=1}^{\infty}a_n \sin(2\pi n f t +\phi_n)
\end{equation}
of the given continuous periodic wave and associate the series (if it converges in the CIE space)
\begin{equation}\label{color_coordinates}
a_0+\sum\limits_{n=1}^{\infty}a_n XYZ(h(n f)),
\end{equation}
but it is to be determined whether (1) this procedure coincides with that used for FM synthesis and (2) the phase $\phi_n$ affects the color quality. These are open questions. In Section~\ref{example} we exemplify computationally the procedure for the FM case. 
 
\section{Categorical enrichment}\label{categorical_enrichment}

Color and timbre, and their relations, can be recast in a categorical framework, where we emphasize the color and timbre transitions, rather than the objects {\em color} and {\em timbre} themselves.

Each topological space $X$ (like the CIE and timbre space) has an associated category whose morphisms are invertible, that is, a \textit{groupoid}. Its objects are the elements of $X$ and its morphisms are homotopy classes of paths in $X$. The composition $[\tau]\circ [\sigma]$ of two classes $[\sigma]:x\longrightarrow y$ and $[\tau]:y\longrightarrow z$ is the class of the concatenation $\sigma\tau$. The identity on $x$ is the class of the associated constant map and the inverse of a path $\sigma$ sends $t\in [0,1]$ to $\sigma(1-t)$. This construction can be generalized to yield higher relations between paths as follows.

\subsection{Induced infinity-groupoids}

\begin{figure}[ht!]
\centering
\includegraphics[width=6cm]{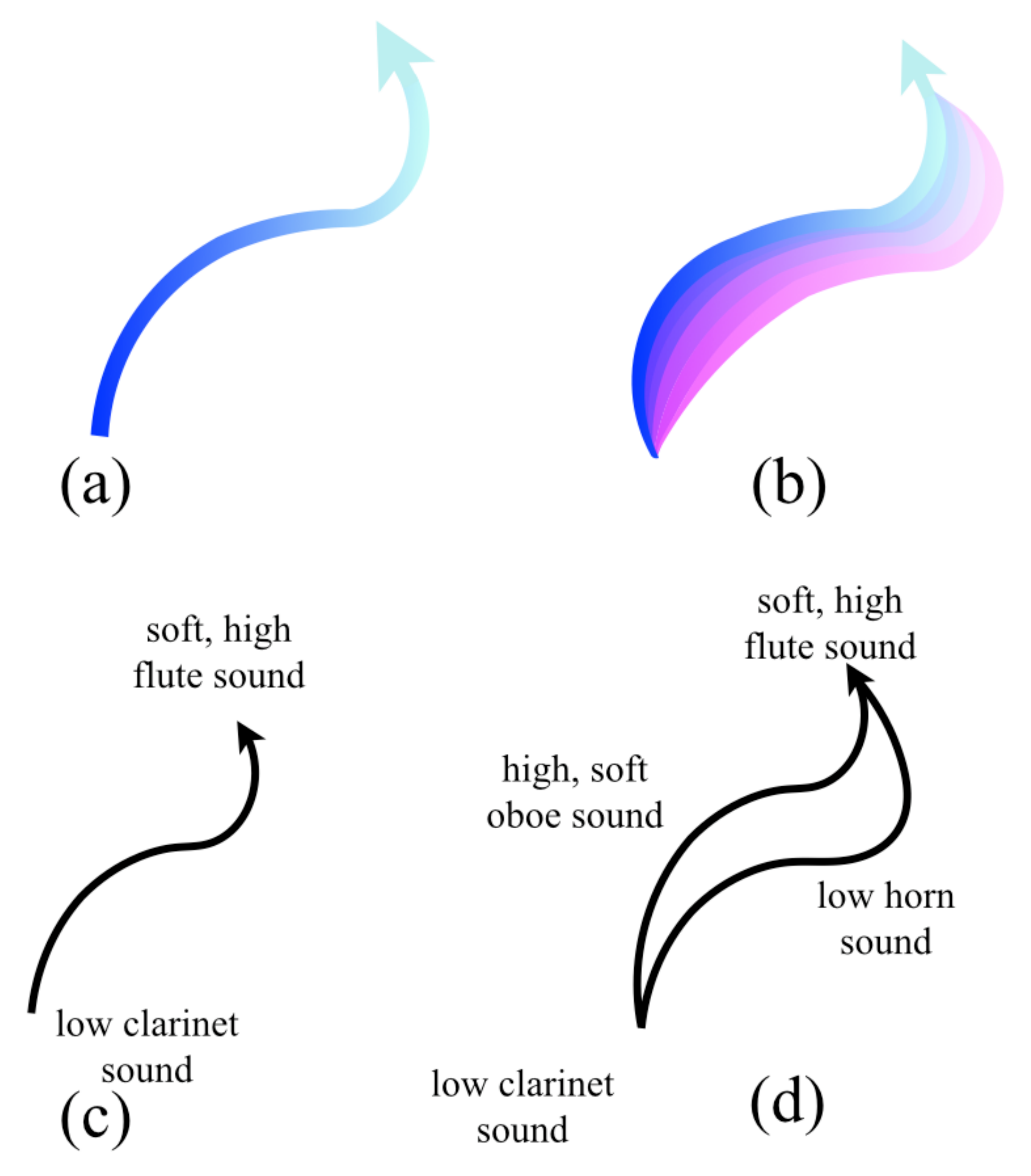}
\caption{(a) A $1$-morphism in the space of colors, a path between two colors, (b) a $2$-morphism in the same space, a band between two color paths, (c) a $1$-morphism in the timbre space, and (d) a $2$-morphism in the same space.} \label{bands}
\end{figure}

Let us consider the \textbf{singular complex} $\s{X}$, which is a \textbf{simplicial set} and an \textbf{$\infty$-groupoid}, under the definitions in Section~\ref{G}.
\textcolor{black}{According Proposition 1.9 and Remark 1.10 from \cite{groth}, $\infty$-categories have $n$-morphisms for each $n\geq0$ and composition of them, which is associative up to homotopy.
}
% This leads to the following sub-definitions.
\textcolor{black}{Thus, \textit{$1$-morphism} of $\s{X}$ is a path in $X$, and a \textit{$2$-morphism} is a homotopy between two paths with the same endpoints.
%However, this is not precise, because there are no Sing$(X)$ 2-morphisms.
} Note that the groupoid of $X$ comes from homotopy classes of $1$-morphisms and hence the concatenation of them is associative up to homotopy equivalence.
On the other hand, a $2$-morphism can be seen as a band of intermediate paths between two given ones that connect the same points. Figure~\ref{bands} shows examples of $1$-morphisms and $2$-morphisms in the cases of the CIE and timbre spaces. 
More generally, we can define $n$-morphisms of the singular complex, which describe the evolution of a single color (timbre), of a path of colors (timbres), of a homotopy of paths, and so on.

We emphasize the need for higher relations and bands. For example, we can map the transition {\em light blue}$\rightarrow${\em dark blue} into the transition {\em light green}$\rightarrow${\em dark green}, creating a band that connects, as different shades, light green with light blue, and dark green with dark blue. \footnote{This relation is not a proper morphism but it can be achieved by pasting two suitable $2$-simplices in $\s{X}$, which is a higher gesture according to \cite{simpgest}. It is also a hypergesture in the sense of \cite{MazzolaGest}.} 
If the initial and final points of the band coincide, we can have the situation described in Figure \ref{bands}, where the dark blue becomes a light blue through different paths: some paths remain in the blue area, while other ones cross the violet area \cite{color_journal}.

\subsection{Induced functors}\label{induced}

Given two topological spaces $X$ and $Y$, which can be the timbre and the CIE color space respectively, and a continuous map $f:X\longrightarrow Y$ there is an induced natural transformation $F:\s{X}\longrightarrow \s{Y}$ that sends a \textit{singular $n$-simplex} $\sigma:\Delta^n\longrightarrow X$ to $f\sigma:\Delta^n\longrightarrow Y$. According to the definition in Section~1.2.7 of \cite{lurie}, which says that a functor between \textbf{infinity-categories} is a natural transformation between the respective simplicial sets, $F$ is a functor from $\s{X}$ to $\s{Y}$.

Note that $F$ coincides with $f$ on objects and sends a $1$-morphism $\sigma$ in $X$ to the path $f\sigma$ in $Y$. 

As any functor between infinity-categories, $F$ preserves the usual categorical structure (up to homotopy), %\textcolor{black}{However, Sing$(X)$ does not have homotopy classes, nor does it have composition. Thus, we should rather think of functors between bigroupoids \cite{color_journal}.} 
in the sense that
\[F([id_x])=[id_{f(x)}]\]
whenever $x \in X$ and 
\[F([\tau]\circ[\sigma])=F([\tau])\circ F([\sigma])\]
whenever $\sigma:x\longrightarrow y$ and $\tau:y\longrightarrow z$ are paths in $X$. More generally, $F$ preserves the compositions of higher morphisms in an appropriate sense, but we omit these technical details. Next, a computational sketch of a functor from timbre to color.

\subsection{A computation of colors from timbres}\label{example}

As an example of associations between a timbre path and a color path, let us consider the progressive enrichment of a simple $440$ Hz sine wave with harmonics, using FM synthesis, and the associated color transition. 

\begin{figure}[ht!]
    \centering
    \begin{minipage}{0.45\textwidth}
        \centering
        \includegraphics[width=\textwidth]{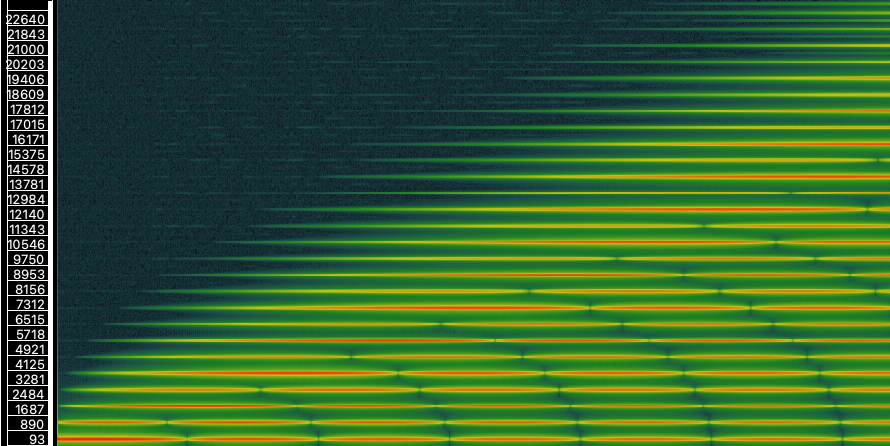} % first figure itself % sound_color_gradient
        \caption{An example of timbre path. The spectrogram is obtained with SonicVisualiser. \textcolor{black}{The darker the color, the closer the sound to silence.}}\label{sound_color_gradient}
    \end{minipage}\hfill
    \begin{minipage}{0.45\textwidth}
        \centering
        \includegraphics[width=0.55\textwidth, angle =-90 ]{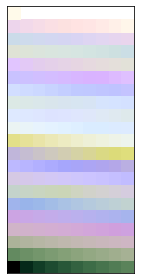} % second figure itself
        \caption{Visual color gradient corresponding to the timbre path of Figure \ref{sound_color_gradient}. \textcolor{black}{Each color corresponds to a value of the modulation index $I$.} % According to the timbre-color correspondence in equations \eqref{bessel1} and \eqref{color_coordinates}.
        }\label{color_gradient}
    \end{minipage}
\end{figure}

More formally, take $f_c=440$ and $f_m=2f_c$. By regarding the modulation index $I$ as a parameter in the interval $[0,20]$, we obtain a continuous\footnote{This defines a continuous map $[0,20]\times \mathbb{R}\longrightarrow \mathbb{R}$, so the exponential transpose $[0,20]\longrightarrow \mathbb{R}^\mathbb{R}$, which is a path, is continuous.} path in the timbre space with parametrization (eq. \ref{bessel1}):
\[ \sin[\omega_c t+I\sin(\omega_m t)]  .\] 
The result is a fluctuation in the brilliance of a sort of clarinet sound since only odd harmonics are present.\footnote{The audio file, in which we identify the increasing modulation index with time (in seconds), can be accessed from the link \url{https://soundcloud.com/maria-mannone/fm-path/s-cFJ7kNrqJjs?}.} Figure~\ref{sound_color_gradient} is the corresponding spectrogram of the timbre path.

To obtain a color path (Figure \ref{color_gradient}) we use the procedure in Section~\ref{map} for each value of $I$, \textcolor{black}{see eq. \eqref{color_coordinates}. For each new value of the index modulation $I$, harmonics vary, reaching a new timbre in Figure \ref{sound_color_gradient}. For each new value of $I$, and thus, for each timbre point reached, there is a color point reached in Figure \ref{color_gradient}. In fact, each color bar represents a color point in the space of colors.} This could mean that we are using the functor induced by any of the continuous maps from timbres to colors (Section~\ref{map}), according to Section~\ref{induced}. In Figure~\ref{color_gradient}, the color squares correspond to the modulation index $I$ values $n/10$ for integers $n$ from $0$ to $200$. There, the modulation index increases from left to right and from top to bottom.
% the possible continuous map from timbres to colors (Section~\ref{map}), according to Section~\ref{induced}.
Then one uses conversion to RGB for screen representation (Section~\ref{color}). 

The Python codes for the FM path and color path \textcolor{black}{are available at \url{https://github.com/medusamedusa/color_gesture}}.

The results agree with Caivano's reflections \cite{caivano}: the closer the sound to a white noise, the closer the color to white light, with additive color mixing. The inverse choice could associate the richness of harmonics (especially in a low-register orchestral range) with a darker color, more like in painting, with subtractive color mixing. In the first case, primary colors are red, blue, and green, and their sum gives white; in the second case, primary colors are red, yellow, and blue, and their sum gives black. In gestural chromo-similarity (Section~\ref{gestsimi}), in analogy with painting we may use the second option (subtractive), see an example in \cite{color_journal}. 

\section{Gestural considerations}\label{gest}

We close this paper with some gestural reflections that may enrich the color and timbre relation theory.

\subsection{From paths to gestures and gestural similarity}\label{gestsimi}

Color and timbre paths (or 1-morphisms) are particular cases of \textbf{gestures}  \cite{gestfun,simpgest,hughes,knot,mathanddance,MazzolaGest,clark}, which are informally diagrams (shaped by a digraph) of paths in a topological space. Continuous maps induce new ones between respective spaces of gestures, as we explain in Section~\ref{section_maps}, so there are correspondences between color gestures and timbre gestures.

We can talk of \textit{gestural similarity} if musical sequences (auditory domain) and simple sketches (visual domain) appear as being produced by the same generator gesture \cite{gestsim}. This possible definition is supported by the hypothesis of a {\em supramodal brain} \cite{rosenblum}. Thus, when gestures in the space of colors and gestures in the space of timbres show perceptive analogies, we can talk of {\em chromo-gestural similarity}.

\subsection{Induced maps between spaces of gestures}\label{section_maps}

Let $\Gamma$ be a digraph. A continuous map $f:X\longrightarrow Y$ induces a new one\footnote{In essence, it is continuous because each component (composition with $f$ for arrows and $f$ for vertices) is.} between \textbf{topological spaces of $\Gamma$-gestures}, namely
\[\Gamma\pitchfork F:\Gamma\pitchfork S_X\rightarrow \Gamma\pitchfork S_Y: \left((c_a)_{a\in A},(x_v)_{v\in V}\right)\mapsto \left((fc_a)_{a\in A},(f(x_v))_{v\in V}\right),\]
\textcolor{black}{where $\Gamma\pitchfork S_X$ ($\Gamma\pitchfork S_Y$) is the \textbf{space of $\Gamma$-gestures in $X$ ($Y$, respectively)}.}

As an example, the Attack-Delay-Sustain-Release (ADSR) envelope of a sound is a gesture shaped by the digraph $\bullet \rightarrow \bullet \rightarrow \bullet \rightarrow \bullet \rightarrow \bullet$ in the amplitude-time space. The envelope has a main role in timbre perception. We can transfer the envelope to the color space by regarding it as an intensity gesture of a single color. 
In fact, this remark may raise new questions regarding color envelopes, and transitions effects from a color to another one.

Color and timbre ramifications, which are gestures, are interesting objects to study and apply to composition. Shaping the orchestral colors, in particular, is a distinctive mark of a composer's style, of a genre, of an epoch. Thus, the proposed ideas can be developed in terms of machine learning as exploited in music information retrieval. Vice versa, a creative interface may be developed starting from the proposed theoretic tools.

Note that the objects and $1$-morphisms of the $\infty$-groupoid $\s{\Gamma\pitchfork S_X}$ are $\Gamma$-gestures in $X$ and paths between gestures, respectively. This groupoid allows one to generalize the idea of timbre paths to transformations between timbres with different ADSR envelopes (loudness profile over time). We may, for example, keep the timbre of a musical instrument while changing its envelope, or keep the envelope and change the timbre, thus performing separate transformations of the envelope and timbre in terms of spectral superposition. As a final abstraction, $\Gamma\pitchfork F$ induces a functor between $\infty$-groupoids $\s{\Gamma\pitchfork S_X}\longrightarrow \s{\Gamma\pitchfork S_Y}$, which would \textcolor{black}{help transfer} envelope transitions between the color and timbre domains.

\section{Conclusion}\label{conc}

The proposed categorical framework could be a way to understand the relation between color and timbre, complementing classical approaches from physics. This framework is based on structural analogies between the perceptual domains of hearing and vision. It is interesting to ask to what extent categorical models could be independent from perception and classical models, taking into account the computational advantages of the latter.  

We also proposed a gestural extension of the categorical framework to capture gestural similarities between the musical and auditory domains.

\textcolor{black}{As a possible, alternative structure to look at, we could consider the Moore paths as 1-cells, in order to have a strictly associative composition, taking homotopies of homotopies for the 2 cells \cite{grandis2}. Given that we are interested in invertible arrows, another suitable structure appears to be the \textbf{bigroupoid} \cite{hardie}, that is, a weakly-invertible \textbf{bicategory}. Concerning the spaces, we could also consider the Euclidean space of colors (as RGB) and the Euclidean space of timbres as defined by Grey \cite{grey}. In a (bi)groupoid, all arrows are invertible. In this way, the points (single colors, single timbres) are 0-cells; the color gestures and timbre gestures are the paths, the 1-cells; the bands (hypergestures in the sense of \cite{MazzolaGest}) are the 2-cells. Path associativity is verified for equivalence classes of homotopies. The model of bigroupoid for color and timbre gestures is discussed in detail in \cite{color_journal}.}

\textcolor{black}{However, we stress the fact that $\infty$-categories simplify the involved axioms and computations in higher category, 2- and bi- categories included.}

This research could lead to signal processing practical implementations, and it could provide a theoretical framework to analyze experiments in the domain of musical timbre.
On the creative side, other possible directions may involve the development of interfaces for composers to manipulate timbres through symbols and/or color references, and for visual artists to do the inverse.

The possibility of translating structures from one domain to another one, provided that some cognitive conditions are verified \textcolor{black}{\cite{gestart}}, can open scenarios also for disability studies, where people with visual impairment can benefit from auditory-accessible interfaces, and people with auditory impairment can benefit from visually-accessible interfaces \cite[pp.~128-129]{caivano}. The reference to gesture and touch regarding intensity, organization, and time distribution of stimuli can inspire even more audacious applications for touch-based interfaces for deaf-blind people.

Thus, a simple question such as ``can we join timbres and colors?'' can open the way to striking applications to improve people's lives.

\section{Glossary}\label{G}

\textbf{Bicategory} \textcolor{black}{In a bicategory, the morphism composition is not associative, but only associative up to an isomorphism. This notion has been introduced by B\'{e}nabou in 1967 \cite{benabou}. The objects are the 0-cells, the morphisms are the 1-cells, and the morphisms between morphisms are the 2-cells.}

\textbf{Bigroupoid} \textcolor{black}{A bigroupoid is a bicategory whose ``2-cells are strictly invertible, and the 1-cells are invertible up to coherent isomorphism” \cite{hardie}.}

\textbf{Compact-open topology} The subbasic opens of the \textit{compact-open topology} on the space of continuous maps $\mathbb{R}^{\mathbb{R}}$ are those of the form \[\{f:\mathbb{R}\longrightarrow \mathbb{R}\text{ continuous }\ |\ f(K)\subseteq U\},\] where $K$ is compact (closed and bounded) in $\mathbb{R}$ and $U$ is open in $\mathbb{R}$. This makes $\mathbb{R}^{\mathbb{R}}$ an exponential in the category of topological spaces. \textcolor{black}{The fact that Top is not Cartesian closed does not imply the non-existence of $\mathbb{R}^\mathbb{R}$.}
\newline

\noindent\textbf{Simplicial category} Denote by $[n]$ the ordered set (ordinal) $\{0,1,\dots,n\}$ for $n\in\mathbb{N}$. The \textit{simplicial category} $\Delta$ has as objects all $[n]$ for $n\in \mathbb{N}$ and as morphisms all order-preserving maps between them. 
\newline

%\noindent \hypertarget{ss}{\textbf{Standard simplex (functor)}}
\noindent\textbf{Standard simplex (functor)}
For each $n\in \mathbb{N}$, we define the \textit{standard $n$-simplex} $\Delta^n$ as the set \[\{(t_1,\dots,t_n)\ | \ 0\leq t_1\leq \dots \leq t_n\leq 1\}.\]
The standard $n$-simplex is a subspace of $\mathbb{R}^n$ and this construction defines a \textit{standard simplex functor} $\Delta^{(-)}$ from the
{\bf simplicial category}
%\hyperlink{sc}{simplicial category}
to the category of topological spaces, which sends an order-preserving map $\alpha:[n]\longrightarrow [m]$ to the appropriate continuous map $\Delta^{\alpha}:\Delta^n\longrightarrow \Delta^m$ sending the $i$th vertex (with $n-i$ zeros) to the $\alpha(i)$th one. \textit{Examples}: $\Delta^{0}$ is a singleton, $\Delta^{1}$ is the interval $[0,1]$ in $\mathbb{R}$; $\Delta^{2}$ is the triangle with vertices $(0,0)$, $(0,1)$, and $(1,1)$ in $\mathbb{R}^2$; and $\Delta^{3}$ is the tetrahedron with vertices $(0,0,0)$, $(0,0,1)$, $(0,1,1)$, and $(1,1,1)$ in $\mathbb{R}^3$.  
\newline

\noindent \hypertarget{sset}{\textbf{Simplicial set}}
Functor from the opposite $\Delta^{op}$ of the {\bf simplicial category} to the category $\mathbf{Set}$ of sets. \textit{Example}: The %\hyperlink{sing}{singular complex}
% \hyperlink{sc}{simplicial category}
{\bf singular complex}
$\s{X}$ of a topological space $X$. 
\newline

%\noindent \hypertarget{sing}{\textbf{Singular complex}}
\noindent\textbf{Singular complex}
The \textit{singular complex of a topological space $X$}, denoted by $\s{X}$, is the {\bf simplicial set} 
%\hyperlink{sset}{simplicial set} 
$\mathbf{Top}(\Delta^{(-)},X)$, where $\Delta^{(-)}$ is the {\bf standard simplex functor}. % \hyperlink{ss}{standard simplex functor}
\textit{Examples}: a $0$-simplex of $\s{X}$ is a point of $X$, a $1$-simplex of $\s{X}$ is a path in $X$.
\newline

%\noindent \hypertarget{infcat}{\textbf{Infinity-category}}\hyperlink{sset}{Simplicial set}
\noindent\textbf{Infinity-category}
A {\bf simplicial set} $S$ is a set such that given $n\in \mathbb{N}$ and $k$ with $0<k< n$, for each subset $\{a_i\ | \ 0\leq i\leq n;\ i\neq k\}$ of $S([n-1])$ satisfying the identities
\[d_i(a_j)=d_{j-1}(a_i)\ \ (i<j;\ i,j\neq k)\],
there is an element $a\in S([n])$ such that $d_i(a)=a_i$ for $i\neq k$. If this property also holds for $k=0$ and $k=n$, then we say that $S$ is an \textbf{$\infty$-groupoid}. \textit{Example}: The {\bf singular complex} of a topological space is an $\infty$-groupoid.
\newline
% \hypertarget{infgr}{\textbf{$\infty$-groupoid}}
% \hyperlink{sing}{singular complex}

%\noindent \hypertarget{gestsp}{\textbf{Topological space of gestures}}
\noindent\textbf{Topological space of gestures}
Let $\Gamma$ be a digraph $(A,V,d_0,d_1)$ and $X$ a topological space. The \textit{space of $\Gamma$-gestures in $X$}, denoted by $\Gamma\pitchfork S_X$ (\textcolor{black}{where $\pitchfork$ stands for transversality}), is the subspace of the product space (compact-open topology on $X^I$)
\[\left(X^I\right)^A\times  X^V\]
consisting of all sequences $\left((c_a)_{a\in A},(x_v)_{v\in V}\right)$ such that $c_a(i)=x_{d_i(a)}$ for $i=0,1$. We say that such a sequence is a \textbf{$\Gamma$-gesture in $X$}.  

\section*{Disclosure statement}% The * makes this section unnumbered.
\vspace{-8pt}
\addcontentsline{toc}{section}{Disclosure statement}% This command makes the unnumbered section also appear in the pdf bookmarks.
No potential conflict of interest was reported by the authors.
\vspace{-8pt}
\section*{Author contributions}
\addcontentsline{toc}{section}{Author contributions}
\vspace{-8pt}
The authors have contributed equally.
\vspace{-8pt}
\nopagebreak[4]

\end{document}